\documentclass[11pt,a4paper]{article}
\usepackage{amsmath,amsthm,enumerate}
\usepackage{amssymb}
\usepackage{amsfonts}
\usepackage{graphicx}
\usepackage{color}
\usepackage{csquotes}
\usepackage{tcolorbox}

\newtheorem{Theorem}{Theorem}[section]
\newtheorem{Example}{Example}[section]
\newtheorem{Proposition}{Proposition}[section]

\newtheorem{Lemma}{Lemma}[section]
\newtheorem{Remark}{Remark}[section]

\newcommand{\sol}{\mathrm{sol\, }}
\def\bc{\begin{center}}
\def\ec{\end{center}}

\def\s2c{\vskip 2cm}
\def\disp{\displaystyle}
\def\cone{\mbox{\rm cone}\,}

\def\gph{\mbox{\rm gph}\,}
\def\epi{\mbox{\rm epi}\,}

\def\dom{\mbox{\rm dom}\,}
\def\bd{\mbox{\rm bd}\,}

\def\R{\mathbb{R}}

\def\vcR{\R \cup \{ +\infty\} }
\def\tR{\Tilde\R^T_+}
\def\cR{\Tilde\R^C_+}
\def\lm{\lambda}

\def\Th{\Theta}
\def\N{\mathbb{N}}

\def\vt{\vartheta}
\def\ox{\bar{x}}

\def\Tilde{\widetilde}

\def\ra{\rangle}
\def\rar{\rightarrow}
\def\la{\langle}

\numberwithin{equation}{section}

\def\bt{\begin{Theorem}}
\def\et{\end{Theorem}}
\def\bl{\begin{Lemma}}
\def\el{\end{Lemma}}
\def\bcor{\begin{Corollary}}
\def\ecor{\end{Corollary}}
\def\supp{\mbox{\rm supp}\,}

\begin{document}

\title{Simple Bilevel Programming and Extensions Part-I: Theory }
\author{Stephan Dempe\footnote{Faculty of Mathematics and Computer Science, TU Bergakademie Freiberg, Germany, e-mail: dempe@tu-freiberg.de. Work of this author has been supported by Deutsche Forschungsgemeinschaft}, Nguyen Dinh\footnote{Department of Mathematics, International University, Vietnam National University-Ho Chi Minh City,
		Ho Chi Minh City, Vietnam, e-mail: ndinh02@gmail.com.  \ \ Work of this author was supported by the project B2019-28-02: Gerneralized scalar and vector Farkas-type results with applications to optimization theory, from Vietnam National University - Ho Chi
Minh city, Vietnam.
}, Joydeep Dutta\footnote{Department of Economic Sciences, Indian
Institute of Technology, Kanpur, India, e-mail: jdutta@iitk.ac.in} and Tanushree Pandit\footnote{Department of Mathematics and Statistics, Indian
Institute of Technology, Kanpur, India, e-mail: tpandit@iitk.ac.in}}
\date{}
\maketitle

\begin{abstract}
 In this paper we begin by discussing the simple bilevel programming problem (SBP) and its extension the simple mathematical programming problem under equilibrium constraints (SMPEC). Here we first define both these problems and study their interrelations. Next we study the various types of necessary and sufficient optimality conditions for the (SMPEC) problems; which occur under various reformulations. The optimality conditions for the (SBP) problem are special cases of the results obtained here when the lower level objective is the gradient of a convex function. Among the various optimality conditions presented here are the  sequential optimality conditions, which do not need any constraint qualification. We also present here a schematic algorithm for the (SMPEC) problem, where the sequential optimality conditions play a key role in the convergence analysis. 

\end{abstract}

\\
\medskip


\section{Introduction and Motivation}\label{intro}
This paper is motivated by the earlier work of the first three authors
on the simple bilevel programming problem \cite{DDD10}. The simple bilevel programming problem consists of minimizing a convex function over the solution set of another convex optimization problem. More formally we can state it as follows. Consider
the following problem
\begin{eqnarray*}
\min f(x), \quad\mbox{subject to}\quad x \in S,
\end{eqnarray*}
where
\begin{eqnarray*}
S = \mbox{argmin} \{h(x) : x \in C \},
\end{eqnarray*}
and $ f $ and $ h $ are real-valued convex functions on
$\mathbb{R}^n $ and $ C $ is a closed convex set. We call it the
simple bilevel programming problem (SBP) since in the original formulation
of a bilevel programming problem there are two variables, one
for the upper-level problem and the other for the lower level
problem while we have only one decision variable, see \cite{dem02}
for details. It is simple to observe that, if we assume that $ h $ is
differentiable, the set $ S $ can be equivalently written as
\begin{eqnarray*}
S = \{ x \in C : 0 \in \nabla h(x) + N_C(x) \},
\end{eqnarray*}
where $ N_C(x) $ is the normal cone to the convex set $ C $ at $ x
$. It is important to note that, unlike the usual bilevel programming
problem which is not in general convex even if the data is convex,
the simple bilevel problem is a convex optimization problem. This
problem was first studied by Solodov \cite{Solodov1} who developed an
algorithm for solving such problems. Solodov \cite{Solodov1} also
showed that many classes of problems can be modeled as a simple
bilevel programming problem. In fact the standard convex
optimization problem can also be modeled as a simple bilevel
programming problem. The problem (SBP) was first theoretically
analyzed in \cite{DDD10} where an approach was developed to generate
necessary and sufficient optimality conditions using very recent techniques from the theory of convex optimization. There were slight discrepancy in the analysis of one result in \cite{DDD10}, which was corrected in the monograph by Dhara and Dutta \cite{Dhara-Dutta}. It is important to note that in \cite{DDD10}, the lower-level problem in (SBP) was allowed to have non-differentiable data. In that case we have $ S = \{ x \in \mathbb{R} : 0 \in \partial h(x) + N_C(x) \}$. In this paper when we consider the (SBP) problem from a theoretical perspective we will always allow the lower-level objective, constraints and the upper-level objective to be non-smooth. We take this opportunity to state two interesting examples of the simple bilevel programming problem.

\begin{Example}
\label{example1.1}  \rm Let us consider the convex problem (CP) in the lower level of (SBP), i.e.
\begin{eqnarray*}
\min h(x), \quad\mbox{subject to} \quad x \in C,
\end{eqnarray*}
where as before, h is a finite-valued convex function on $\mathbb{R}^n$. It is meaningful to ask the question if, given an $x \in C$, can we estimate $d(x,S)$, where $S$ denotes the solution set of (CP), or at least provide error bounds? When $ h $ is strongly convex one can devise global error bounds for $ d (x,S)$ using the machinery of gap functions by viewing the problem (CP) as a variational inequality problem. See for example Fukushima \cite{Fuku}. However if $h$ is just convex and need not be strongly convex, then the problem is far from being resolved. But we can still approach the problem by posing it as a simple bilevel problem. Using the (SBP) formulation we can numerically estimate $ d(x,S)$ for a given $x$, rather than finding global error bounds. However the estimation of $ d(x,S) $ is important in practice. For a given $x \in C$ in order to estimate $d(x,S)$ we need to solve the problem
\begin{eqnarray*}
\min \frac{1}{2}\|y-x\|^2, \quad\mbox{subject to} ~ y \in S; \quad \mbox{where} ~ S= \mbox{argmin}\{ h(x): x \in C \}.
\end{eqnarray*}
Thus we have posed our problem as a simple bilevel problem. In fact it might appear strange that in order to compute $d(x,S)$ we need to have a knowledge of $S$. \textit{However our algorithm will show that such a requirement is not necessary.} Thus the simple bilevel programming approach is effective to compute $d(x,S)$.\\\phantom{b} \hfill $\Box$
\end{Example}

\begin{Example}
\label{example1.2} \rm Another interesting case where at least the upper level has non-smooth data is given as follows. Consider the problem
\begin{eqnarray*}
\min \|x\|_1 \quad\mbox{subject to} \quad Ax=b
\end{eqnarray*}
where $A$ is a $m \times n$ matrix (may be of full row rank) and $b \in \mathbb{R}^m$. This is often called the basis-pursuit problem which provides sparse solution to the system $Ax=b$. This can be posed as the simple bilevel problem
\begin{eqnarray*}
\min \|x\|_1 \quad\mbox{subject to} \quad x\in \mbox{arg}\min\limits_{x \in \mathbb{R}^n}\|Ax-b\|^2.
\end{eqnarray*}
This is clearly a simple bilevel programming problem. $\hfill$ $\Box$
\end{Example}

If we consider the problem (SBP) in which the lower level problem has a smooth objective function then it can be easily generalized as follows. Consider the following problem which we will call (SMPEC);
\begin{eqnarray*}
\min f(x), \quad\mbox{subject to}\quad x \in S,
\end{eqnarray*}
where $f$ is a continuous convex function and $ S=\sol(VI(F,C)) $ is the solution set
of the variational inequality $VI(F, C) $, where we seek to find $ x
\in C $ such that
\begin{eqnarray*}
\langle F(x), y - x \rangle \ge 0,  \quad \forall y \in C,
\end{eqnarray*}
where $ F : \mathbb{R}^n \rightarrow \mathbb{R}^n$ is a continuous, monotone map and $ C $ is a closed convex set in $ \R^n $. In the sequel, we always assume that $S$ is a
non-empty  set. If $ C $ is a compact convex set then $ S $
is always non-empty. It is important to note that if $ F = \nabla h $, where $ h $ is a convex function, then (SMPEC) reduces to (SBP) with a smooth lower level objective. Thus (SMPEC) is a generalization of (SBP) when the lower level objective of (SBP) is smooth. It is crucial to note that if the lower level problem of (SBP) has non-smooth objective function then the SMPEC problem mentioned above cannot be considered as a generalization of SBP.  \\
We would like to note that even if the lower-level objective function is not smooth, then we call still have a generalization of the (SBP). This is given as follows. Consider the following problem
\begin{eqnarray*}
\min f(x) \quad \mbox{ subject to } ~ 0 \in T(x)+ N_C(x),
\end{eqnarray*}
where $ T: \mathbb{R}^n \rightrightarrows \mathbb{R}^n $ is a maximal monotone map. This problem is referred to as the simple mathematical problem with GVI constraints which we can denote as (SMPGVI). If $ T \equiv \partial f $ where $ \partial f $ is the subdifferential of a convex function, then (SMPGVI) reduces to the problem (SBP) with a nonsmooth lower-level objective. We will take up the study (SMPGVI) in a separate paper.\\
Further we would like to draw the attention of the reader to section 3.2 where it is shown that SMPEC problem presented here can be reformulated as a simple bilevel problem with non-smooth lower level objective function.\\
\noindent We call the above (MPEC) problem the simple MPEC problem for reasons 
similar to that of the simple bilevel problem. The MPEC or {\it
Mathematical Program with  Equilibrium Constraints} is studied in
detail in Luo, Pang and Ralph \cite{LuPaRa96}. In their setting the
variational inequality problem at the lower-level is a parametric one,
i.e.,  the MPEC problem like the bilevel problem has
two variables.

The first question is whether the problem (SMPEC) is a
convex optimization problem. The answer is yes if $F$ is a
continuous and monotone map, that is,
 for any
$ x, y \in \mathbb{R}^n $
\begin{eqnarray*}
\langle F(y) - F(x), y - x \rangle \ge 0,
\end{eqnarray*}
(see for example Facchinei and Pang \cite{fp-2003}). Then it is meaningful to
talk about necessary and sufficient optimality conditions for the problem (SMPEC).  The details of this will be discussed in this article.

Here we provide a motivating example of an (SMPEC) problem.
\begin{Example} \rm
Consider the following pair of primal - dual linear programming problems
\begin{eqnarray*}
\mbox{(LP)} &: &\min \langle c,x \rangle,\quad \mbox{subject to}\quad Ax \geq b, x \geq 0\\
\mbox{(DP)} &: & \max \langle b,y \rangle,\quad \mbox{subject to}\quad A^Ty \leq c, y \geq 0,
\end{eqnarray*}
where $A$ is a $m \times n$ matrix, $c \in \mathbb{R}^n$ and $b \in \mathbb{R}^m$.

Suppose that  we consider the problem of finding the primal-dual solution of minimum norm. Consider the $VI(F(x,y), \mathbb{R}_+^n \times \mathbb{R}_+ ^m) $ given as

\begin{eqnarray}\label{1.1}
F(x,y)=
  \begin{bmatrix}
    0 & -A^T \\
    A & 0
  \end{bmatrix}
  \begin{bmatrix}
    x\\
    y
  \end{bmatrix}
  +
  \begin{bmatrix}
    c\\
    -b
  \end{bmatrix}.
\end{eqnarray}
Let us set
\[
M=
\begin{bmatrix}
0 & - A^T \\
A & 0
\end{bmatrix}.
\]
Then $M$ is a skew-symmetric matrix and hence positive semidefinite:
\[\langle (x,y), M(x,y) \rangle =0 \textrm{ for all }   (x,y) \in \mathbb{R}^n \times \mathbb{R}^m.\]
Let $S(A,b,c)$ denote the primal-dual solution set of (LP) and (DP) and let $\sol(VI(F(x,y), \mathbb{R}_+^n \times \mathbb{R}_+^m))$ be the solution set of $VI(F(x,y), \mathbb{R}_+^n \times \mathbb{R}_+^m) $, where $F$ is given by (\ref{1.1}). Then it can be shown that
\begin{eqnarray*}
S(A,b,c) = \mbox{sol}(VI(F(x,y), \mathbb{R}_+^n \times \mathbb{R}_+^m)).
\end{eqnarray*}
For a proof see for example Borwein and Dutta  \cite{lon-jd-2015} or Borwein and Lewis \cite{jon-lewis-2006}. Thus the problem that we had posed above can be written as
\begin{eqnarray*}
\min \|(x,y)\|^2  \quad \mbox{subject to} \quad (x,y) \in \mbox{sol}(VI(F(x,y), \mathbb{R}_+^n \times \mathbb{R}_+^m).
\end{eqnarray*}
This is of course an (SMPEC).

The interesting feature of this example is that, though $ VI(F(x,y), \mathbb{R}_+^n \times \mathbb{R}_+^m) $ is an affine variational inequality, the function $F(x,y)$  cannot be obtained as the gradient of any convex function. 
Thus $VI(F(x,y), \mathbb{R}_+^n \times \mathbb{R}_+^m)$ is not the necessary optimality condition of a convex optimization problem.
\hfill $\Box$
\end{Example}

\noindent As mentioned above we have been motivated to study
the simple MPEC problems as a generalization of the simple bilevel
programming problem. However we would also like to note that the
simple MPEC problem can be viewed as a special case of the {\it
Mathematical Programming Problem with a Generalized Equation
Constraint} called (MPGE) for short. This was studied in Ko\v{c}vara
and Outrata \cite{jiri2004}. The (MPGE) problem is stated as follows
\begin{eqnarray*}
\min f(x)\quad\mbox{subject to}\quad 0 \in F(x) + Q(x), x \in C
\end{eqnarray*}
\noindent where $ F : \mathbb{R}^n \rightarrow \mathbb{R}^n $ and $
Q : \mathbb{R}^n \rightrightarrows \mathbb{R}^n $ is a set-valued
map. When $ Q = N_C $ the above problem is nothing but the simple
MPEC problem. Thus in this article we are interested in studying a
particular class of MPGE problems. However it will be important to
note that simple MPEC problem (SMPEC) is also an important one from
the point of view of applications. This has been recently
demonstrated in Facchinei, Pang, Scutari and Lampariello
\cite{fpang2011}. To be more precise Facchinei et al \cite{fpang2011} considered what they termed as variational inequality constrained-hemivariational inequality problem. Consider a continuous and monotone mapping $F: \mathbb{R}^n \rightarrow \mathbb{R}^n$,  a closed and convex set $C$,  a convex function $\phi : \mathbb{R}^n \rightarrow \mathbb{R}$ and a continuous function $H: \mathbb{R}^n \rightarrow \mathbb{R}^n$. The variational inequality constrained-hemivariational inequality problem, denoted as VI-C HVI$(F,C,H,\phi)$, seeks to find $x \in \sol(VI(F,C))$ such that
\begin{eqnarray*}
\langle H(x),y-x \rangle + \phi(y) - \phi(x) \ge 0, \quad \forall y \in \mbox{sol}(VI(F,C)).
\end{eqnarray*}
Note that when $H \equiv 0$ , then VI-C HVI$(F,C,0,\phi)$ coincides with the (SMPEC) problem with $\phi$ as the upper-level objective function.

The article contributes to two different aspects of the (SMPEC) problem. These two aspects are namely
\begin{tcolorbox}
\begin{enumerate}
\item Necessary and sufficient optimality conditions for (SMPEC) problems,
\item A schematic algorithm for the (SMPEC) problem.
\end{enumerate}
\end{tcolorbox}
We will show in this article that a numerical scheme for solving the (SMPEC) problem can be envisaged in connection to the sequential optimality conditions developed in this article. It was during the revision of the first draft of this paper that we took a much deeper look into the role of the optimality conditions in developing numerical schemes for (SMPEC) problems as suggested by the referees. The sequential optimality conditions have played a major role in this and the numerical scheme that emerged for the (SMPEC) problem leads us to solve an (SBP) problem with non-smooth lower level objective function. 
These issues are discussed in Part-II \cite{part2} which deals with algorithms for non-smooth (SBP) and (SMPEC) problems that we have obtained here. Let us also note that our approach to the optimality
conditions for the simple MPEC problem is quiet different from the
one given in Ko\v{c}vara and Outrata \cite{jiri2004} for the MPGE
problem. First of all we consider that the function $ F $ in the
variational inequality $ VI(F, C) $ is a monotone map. We derive
optimality conditions using several different approaches. In the
first approach we view the problem $VI(F, C) $ as a system of
semi-infinite convex inequalities and then use recent techniques
from the literature of convex semi-infinite problems in order to
derive necessary and sufficient optimality conditions for the case
when $ f $ is convex. Then we reformulate the simple MPEC problem as
single objective convex optimization problem using a convex gap
function associated with $ VI(F, C) $ and use a weaker version of
the basic constraint qualification (see \cite{rocw}) to derive the
necessary and sufficient optimality conditions. It is important to note that all the optimality conditions presented in this article are of the necessary and sufficient type.  However,  this may not always be the case owing to the bilevel structure of the problem (see, e.g.,   \cite[Theorem 1]{DDD10}). In fact the bilevel structure allows us to formulate the simple bilevel problem or simple MPEC problem in several equivalent ways. One of them may lead to an optimality condition which is necessary but not sufficient. This was demonstrated in Dempe, Dinh and Dutta \cite{DDD10} for the case of the simple bilevel problem. In this article we want to focus on those formulation for the simple MPEC problem which leads to necessary and sufficient optimality conditions.

The article is planned as follows. In Section 2 we present some basic definitions and facts which we will use throughout the paper. In Section 3 we present  optimality conditions for a simple MPEC problem. This section is further divided into three subsections. In the first one we discuss how to reformulate the (SMPEC) problem as a semi-infinite programming problem which then leads to a type of necessary and sufficient optimality conditions under some closedness-type qualification condition. In the second subsection, we consider a reformulation of the (SMPEC) problem as a single-level problem using the dual gap function of the monotone variational inequality at the lower level. A necessary and sufficient condition is developed using a very weak constraint qualification called the weak-basic constraint qualification. The last subsection of Section 3 develops an optimality condition under calmness condition. We have not been able to find any apparent relation between the calmness condition and the weak basic constraint qualification condition. In Section 4, we deal with sequential optimality conditions for (SMPEC) and also a schematic algorithm for the (SMPEC) problem based on a sequential optimality condition. Note that in deriving these sequential optimality conditions, no constraint qualification is required. In section 4 we divide our study into two cases. The first case does not consider any explicit representation of the convex set $C$ while in the second case $C$ is represented through convex inequality constraints. After we have presented the sequentially optimality conditions we present a numerical scheme for the (SMPEC) problem and the sequential optimality condition is satisfied by the iterates generated in the scheme.
\section{ Preliminaries and Convex Semi-Infinite Programs} \label{sec22}

\subsection{ Preliminaries from Convex Analysis}


The following notations in convex analysis will be used. The {\it
indicator function} $\delta _{D}$ of a set $D$ in $\R^n$ is defined
as $\delta _{D}(x)=0$ if $x\in D$ and $\delta _{D}(x)=+\infty $ if
$x\notin D$.

 Let $f: \R^n \rightarrow
\mathbb{R}\cup \{+\infty \}$ be a proper  function. This means that the set dom$f$, given by $\mbox{dom}f= \{x \in \mathbb{R}^n : f(x) < + \infty\}$ is nonempty. The {\it epigraph} of $f$, epi$f$, is defined by
\begin{equation*}
\mathrm{epi} f=\{(x,r)\in \R^n \times \mathbb{R}\ :\  x\in \mathrm{dom}%
\,{f},\,f(x)\leq r\}.
\end{equation*}%
A function $f: \mathbb{R}^n \rightarrow \mathbb{R} \cup \{+ \infty\}$   is convex iff epi$f$ is convex. Further  it is lower semicontinuous  iff epi$f$ is closed. When $f$ is convex and lower semicontinuous (lsc), the {\it
subdifferential} of $f$,\thinspace\ $\partial f:\R^n
\rightrightarrows \R^n$,  is defined as
\begin{equation*}
\partial {f}(x)=\{v\in \R^n\; :\;  f(y)\geq f(x)+\langle v,y-x\rangle,\,\forall \,y\in
\R^n \},
\end{equation*}%
when $f(x) < + \infty $ and $\partial f (x) = \emptyset $ if $f(x) = + \infty$. Further  if $x \in \mbox{int(dom}f)$, then $\partial f(x) \neq \emptyset$.
It is a well known fact that $\bar{x}$ is a minimizer of $f$ on $\mathbb{R}^n$ if and only if $0 \in \partial f(\bar{x})$. Moreover, if $f$ is a finite-valued convex function then $\bar{x}$ is a minimizer of $f$ over a convex set $C$, if and only if $0 \in \partial f (\bar{x}) + N_C(\bar{x})$, where $N_C(\bar{x})$ is the usual normal cone of convex analysis and $N_C(\bar{x}) = \partial \delta _C (\bar{x})$.

The $\varepsilon $-{\it subdifferential} of $f$, $\partial
_{\varepsilon }f: \R^n \rightrightarrows \R^n$,  is defined as
\begin{equation*}
\partial _{\varepsilon }{f}(x)=\{v\in \R^n\; :\;  f(y)\geq
f(x)+\langle v,y-x\rangle - \varepsilon ,\,\forall \,y\in \R^n \}.
\end{equation*}%
 The {\it
conjugate} function of $f,\ f^{\ast }: \R^n \rightarrow
\mathbb{R}\cup \{+\infty \}$, is defined by
\begin{equation*}
f^{\ast }(v)=\sup \{\langle v,x\rangle-f(x)\ :\  x\in \ \mathrm{dom}\,f\}.
\end{equation*}%
It follows easily from the definitions of $\mathrm{epi}\ f^{\ast }$
of a proper convex function $f$ and the $\varepsilon
$-subdifferential of $f$ that if $a\in \mathrm{dom}f$, then
\begin{equation}\label{epifstar}
\mathrm{epi}\ f^{\ast }=\bigcup_{\epsilon \geq 0}\{(v,\langle v,a\rangle +\epsilon
-f(a))\ :\ \ v\in \partial _{\epsilon }f(a)\}.
\end{equation}%
 For details see \cite{jld}. 
The previous notions are defined by the same way for the infinite
dimensional spaces (see \cite{Zalinescu}).

\begin{Proposition}
\label{proposition2.1} (Brondsted-Rockafellar Theorem \cite{tl}) Let
$X$ be a Hausdorff topological vector space and $a \in X$ and let  $f:X\rightarrow
\mathbb{R}\cup \{+\infty \}$ be a proper, lsc, and convex function.
Then for any real number $\epsilon \geq 0$ and
any $u\in \partial _{\epsilon }f(a)$ there exist $x_{\epsilon }\in X$ and $%
u_{\epsilon }\in \partial f(x_{\epsilon })$ such that
\begin{equation*}
\Vert x_{\epsilon }-a\Vert \leq \sqrt{\epsilon },\Vert u_{\epsilon
}-u\Vert \leq \sqrt{\epsilon }\mbox{ and }|f(x_{\epsilon
})-\langle u_{\epsilon },x_{\epsilon }-a\rangle-f(a)|\leq 2\epsilon .
\end{equation*}
\end{Proposition}
\subsection{Optimality conditions for Convex Semi-Infinite
Programs}

In this section we consider a general class of {\em convex
semi-infinite programs}  (SIP)   of the form:
\begin{eqnarray*}\label{2.1}
{\rm (SIP)} \ \ \left\{\begin{array}{ll}
\mbox{minimize }\;\vt(x)\\
\mbox{ subject to}\\
\ \ \ \ \ \vt_t(x)\le 0,\;\;\forall\,t\in T, \\
\ \ \ \ \ x\in\Th ,
\end{array}\right.
\end{eqnarray*}
where   $T$ is a (possibly infinite) index
set,  $\Th\subset \R^n$ is a closed convex subset,   and
 $\vt\colon \R^n\to\vcR$, 
 and $\vt_t\colon
\R^n \to\vcR$ are proper, lsc, convex functions.

The problem (SIP) is called a {\em convex semi-infinite problem} and
was examined in details in \cite{GL98}. The case where the decision
variables are  in  infinite dimensional space was widely developed
in the recent years (see  \cite{DMN2} and references
 therein).

The set of {\em feasible solutions} to (SIP) will be denoted by
$\Xi$, i.e.,
\begin{eqnarray*}
\Xi:=\Th\cap\big\{x\in X \ :\ \vt_t(x)\le 0\;\mbox{ for all }\;t\in
T\big\}.
\end{eqnarray*}
Further, let $\R^T$ be the {\em product space} of
$\lambda=(\lambda_t\; :\; t\in T)$ with $\lambda_t\in\R$ for all
$t\in T$. Let $\widetilde\R^T$ be collection of $\lambda\in\R^T$ such
that $\lambda_t\ne 0$ for {\em finitely many} $t\in T$, and let
$\tR$ be the {\em positive cone} in $\widetilde\R^T$ defined by
\begin{eqnarray*}
\tR:=\big\{\lambda\in\widetilde\R^T\; :\;\lambda_t\ge 0\;\mbox{ for all
}\;t\in T\big\}.
\end{eqnarray*}
Observe that, given $u\in\R^T$ and $\lambda\in\widetilde\R^T$ and
denoting $\supp\lambda:=\{t\in T \; :\;\lambda_t\ne 0\}$, we have
\begin{eqnarray*}
\lambda u:=\disp\sum_{t\in T}\lambda_t u_t=\sum_{t\in{
 {\rm
supp}}\,\lambda}\lambda_t u_t.
\end{eqnarray*}

The following  qualification condition plays a crucial role in
deriving necessary optimality conditions  for (SIP).

{\it Closedness qualification condition.} We say that the problem
(SIP) satisfies the {\em closedness qualification condition}, (CQC)
in brief,  if the set
\begin{eqnarray*}
\epi\vt^*+\disp\cone\Big\{\bigcup_{t\in
T}\epi\vt^*_t\Big\}+\epi\delta_{\Th}^*
\end{eqnarray*}
is  closed in the space $\R^{n+1}$.

\begin{Remark} It is worth noting that if the  cost function $\vt$
of (SIP) is continuous at some point of the feasible set $\Xi$  or
if the conical set cone $(\dom\vt-\Xi)$ is a closed subspace of $X$,
then the (CQC) requirement  holds provided that the set
\begin{eqnarray*}
\disp\cone\Big\{\bigcup_{t\in T}\epi\vt^*_t\Big\}+\epi\delta_{\Th}^*
\end{eqnarray*}
is  closed in $\R^{n+1}$ (see \cite{DGLS},\cite{DMN2} for more details).
\end{Remark}
\begin{Remark}
Note also that the dual qualification conditions of the (CQC) type
have been introduced and broadly used in \cite{BJ}, \cite{DGLS}, \cite{DMN1}, \cite{DMN2},
\cite{DVN} and other publications of these authors for deriving duality
results, stability and optimality conditions for various constrained
problems of convex and DC programming. It has been proved in the
mentioned papers (see e.g., \cite{DMN1}) that the qualification
conditions of the (CQC) type strictly improved constraint
qualifications of the nonempty interior and relative interior types
for problems considered therein.
\end{Remark}
\noindent A necessary and sufficient optimality condition for (SIP) is given
in the next theorem (see \cite{DMN1}, \cite{DMN2}) which is the key tool for
deriving one of the main results in the next section.

\begin{Theorem} \cite{DMN2}\label{thm21}
Let $\ox\in\Xi$ with $\vt(\ox)<\infty$, and let the qualification
condition (CQC) hold for (SIP).  Then $\ox$ is optimal to (SIP) if
and only if there is $\lm\in\tR$ such that the following generalized
Karush-Kuhn-Tucker $(KKT)$ condition holds:
\begin{equation}\label{21}
0\in \partial \vt(\ox)+\sum_{t\in T}\lambda _{t}\partial
\vt_{t}(\ox)+N_{\Th}(\ox)\ \text{ and}\;\ \lambda
_{t}\vt_{t}(\ox)=0,\;\forall t\in T .
\end{equation}
\end{Theorem}

It is worth observing that for the problem (SIP), active constraints
with positive multipliers
  corresponding to one
of its solutions remain active to all other solutions (see
\cite{SD}). Based on this fact, if one solution of the problem (SIP)
is known, the set of all solutions of (SIP) can be found by a simple
formula as shown in the next theorem. Let $D$ be the solution set of
(SIP) and assume that $D \not=\emptyset$.

\begin{Theorem} \cite{SD}\label{thm22}
For the problem (SIP), suppose that $\bar x$ is a solution of (SIP)
(i.e., $\bar x \in D$) and the (KKT) condition (\ref{21}) holds for
$\bar x$ with the associated Lagrange multiplier $\lm\in\tR$. Then
$D=D_1=\overline{D}_1$, where
\begin{equation*}
\begin{array}{ll}
D_1:=\{x \in \Th\; \ :\  \; &\exists u \in \partial \vt(\bar x)\cap
\partial \vt(x), \la u, x- \bar x \ra=0,\; \vt_t(x)=0, \forall t \in {\rm
supp} \lambda
 , \quad \quad \\
    &\vt_t(x)\le 0, \forall t \not\in  {\rm supp} \lambda\},\\
\overline{D}_1:=\{x \in \Th\; \ :\  \; &\exists u \in \partial
\vt(x), \la u, x- \bar x\ra=0,\;
\vt_t(x)=0, \forall t \in {\rm supp} \lambda , \quad \quad \\
    &\vt_t(x)\le 0, \forall t \not\in  {\rm supp} \lambda\}.
\end{array}
\end{equation*}
\end{Theorem}

\section{Optimality conditions for the Simple MPEC Problem}\label{sect3}

\subsection{SIP approach with a closedness qualification condition}

In this section we show that the simple MPEC problem (SMPEC) can be
converted to a convex semi-infinite  programming problem (SIP). The
results from the theory of convex (SIP) are then applied to get
corresponding ones for (SMPEC). Moreover, at the end of this
subsection, we present a characterization of the solution set of the
(SMPEC) problem. Concretely, we show that whenever one solution
$\bar x$ of (SMPEC) problem is known, based on the Lagrange
multipliers associated to $\bar x$,  an explicit and simple formula
is given that permits one to compute all the other solutions of the
problem in consideration.

{\em  Simple MPEC as a Convex SIP.}   To show this we first need to
show that $ S $ is convex. Using Theorem 2.3.5 in Facchinei and Pang
\cite{fp-2003} one can show that
\begin{equation}\label{31}
S = \bigcap\limits_{y \in C} \{ x \in C\; :\; \langle F(y), y - x
\rangle \ge 0 \}.
\end{equation}
\noindent It is clear that $ S $ is a convex set and thus (SMPEC) is
indeed a convex optimization problem.

This formulation is quite simple by noting the structure of $ S $
given in (\ref{31}). If for  a given and fixed $ y \in C $ we set
\begin{eqnarray*}
\psi_y (x) = \langle F(y), x - y \rangle,
\end{eqnarray*}
\noindent then it is simple to observe that (SMPEC) can be
alternatively formulated as the following  convex SIP
\begin{eqnarray}\label{32}
\left\{\begin{array}{ll}
\mbox{minimize }\; f(x) \; \\
\mbox{ subject to}\\
\ \ \ \ \ \psi_y(x)\le 0,\;\textrm{ for all }y\in C,\;\\
\ \ \ \ \ x\in C.
\end{array}\right.
\end{eqnarray}
 Necessary and sufficient
conditions for (SMPEC) are given in the following theorem.

\begin{Theorem}\label{thm31} Let $\bar x$ be a solution of the variational inequality $VI(F,C)$.
Assume that $F$  is continuous and monotone and the set
\begin{eqnarray}\label{33}
  \ \ \ \ \ \  \cone \bigcup_{y \in C}\{ F(y)\}
\times [\la F(y), y \ra, + \infty)   + \epi \delta_C^\ast
\end{eqnarray}
is  closed in  $\mathbb{R}^{n+1}$. Then $\bar x $ is a solution of
(SMPEC) if and only if there exist some $k\in \mathbb{N}$,
$\lambda_1, \lambda_2, \cdots, \lambda_k > 0 $, $y_1, y_2, \cdots,
y_k \in C$ such that
\begin{eqnarray}
&& 0 \in \partial f(\bar x) + \sum_{i = 1}^{k} \lambda_i F(y_i) +
N_C(\bar x) \ {\rm and} \label{34}\\
&&  \langle F(y_i), \bar x - y_i \rangle = 0, \ {\rm for \ all }\  i
= 1, 2, \cdots , k.\label{35}
\end{eqnarray}
\end{Theorem}
{\bf Proof.} We first note that the problem (SMPEC) can be
re-established as (\ref{32}), which is of the general model (SIP),
where the space of decision variable is $\R^n$,  where $f$ and
$\psi_y$ play the role of $\vt$ and $\vt_t$, respectively, and the
set $C$ in the set constraint  plays the role of the index set as
well.

On the other hand,  by the construction of $\psi_y$, for each $y \in
C$, we have
\[ \psi_y^\ast (x^*) =\left\{\begin{array}{ll} \la F(y),
y \ra \ \ \
&\mathrm{if }\ \ x^* = F(y) \\
+ \infty &\mathrm{ otherwise }.
\end{array} \right.
\]
Consequently,
\[
\epi \psi_y^* = \{ F(y)\} \times [\la F(y), y\ra , + \infty), \] and
hence,
\[
\bigcup_{y \in C} \epi \psi_y^\ast  =\bigcup_{y \in C}\{ F(y)\}
\times [\la F(y), y \ra, + \infty).\]
 Consequently, the assumption
(\ref{33}) assures that the set $\cone \bigcup_{y \in C} \epi
\psi_y^\ast + \epi \delta_C^\ast$ is closed in $\mathbb{R}^{n+1}$,
which, together with the fact that $f$ is continuous (see Remark
2.1), shows that the requirement for the qualification condition (CQC)
for the problem (\ref{32}) holds.

As  $\bar x$ is
a feasible solution of the problem (\ref{32}),
it now follows from Theorem \ref{thm21} that $\bar x$ is a solution
of (SMPEC) (the same, solution of (\ref{32})) if and only if there
exists $\lambda = (\lambda_y)_{y \in C} \in \cR$ such that
\begin{eqnarray}
&& 0 \in \partial f(\bar x) + \sum\limits_{y \in {\rm supp} \lambda}
\lambda_y F(y) +
N_C(\bar x) \ {\rm and} \label{36}\\
&& \lambda_y \langle F(y), \bar x - y \rangle = 0, \ {\rm for \ all
}\  y \in C.\label{37}
\end{eqnarray}
Note that $\supp \lambda$ is finite, we can assume that $\supp
\lambda = \{ y_1, y_2, \cdots, y_k\}$, for some integer $k$. By
setting $\lambda_{y_i} = \lambda_i$, $ i = 1, 2, \cdots, k$,  the
last two expressions, namely, (\ref{36}) and (\ref{37}), become
\begin{eqnarray*}
&& 0 \in \partial f(\bar x) + \sum_{i = 1}^{k} \lambda_i F(y_i) +
N_C(\bar x) \ {\rm and} \\
&&  \langle F(y_i), \bar x - y_i \rangle = 0, \ {\rm for \ all }\  i
= 1, 2, \cdots , k,
\end{eqnarray*}
which is desired. \hfill $\Box$

\medskip
We now apply Theorem \ref{thm22} to characterize the solution set
$S_{M}$ of problem (SMPEC).  It is shown, in particular,
in  the case where $f$ is differentiable and  the set of all
solutions of $VI(F,C)$), sol($VI(F, C)$),   is known,  $S_M$ can be found simply as the
intersection of sol($VI(F, C)$), $C$, and the set of solutions of a finite system of
linear equations.

\begin{Theorem}\label{thm32} Assume that $F$  is continuous and monotone and that  $\bar x$ be a solution of the
problem (SMPEC) such that the  conditions (\ref{34})-(\ref{35}) hold
for $\bar x$ with $k\in \mathbb{N}$, $\lambda_1, \lambda_2, \cdots,
\lambda_k
> 0 $, $y_1, y_2, \cdots, y_k \in C$. Then
\begin{eqnarray*}
{S}_M &=& \Big\{ x\in C \; :\;  \exists u \in \partial f(x ), \la u
, x - \bar x\ra = 0,\Big.\\ &&\Big.\la F(y_i), x - y_i\ra = 0, \ i = 1, 2,
\cdots, k,
 \la F(y), x - y\ra \leq 0, \ \forall y \in C  \Big\}.
\end{eqnarray*}
In particular, if $f$ is differentiable and ${\rm sol}(VI(F, C))$ is known, the
solution set $S_M$ can be characterized simply as:
\[
S_M = \Big\{ x \in S \; :\;  \la \nabla f(x), x - \bar x \ra = 0, \
\la F(y_i), x - y_i\ra = 0, \ i = 1, 2, \cdots, k \Big\} .
\]
\end{Theorem}
{\bf Proof.} This is a direct consequence of Theorem \ref{thm22} and
Theorem \ref{thm31}. \hfill $\Box$

\subsection{The use of a  gap function and the weak basic constraint qualification}

{\it Reconstruction of (SMPEC) using a gap function.} We shall show
that under certain assumptions the simple MPEC problem (SMPEC) can
be posed as a simple bilevel problem and thus the known necessary
and sufficient optimality conditions for the simple bilevel problem
can be used to deduce new optimality conditions for the simple MPEC
problem. Our approach is to reformulate the simple MPEC problem
(SMPEC) into a simple bilevel  problem that rests on our use of a
gap function for $ VI(F, C) $. More specifically we will use the
dual gap function $ g_D $ which has the additional advantage of
always being convex. The dual gap function is given as
\begin{eqnarray*}
g_D (x) = \sup_{y \in C} \langle F(y), x - y \rangle.
\end{eqnarray*}
It is not difficult to see that $ g_D $ is convex and is finite if $
C$ is compact in addition to being convex. Moreover using
 Proposition 2.3.15 from \cite{fp-2003} we can
show that $ g_D $ is a gap function for $ VI(F, C) $ when $ F $ is
continuous and monotone. Since $ g_D $ becomes a gap function it is
not difficult to see that
\begin{eqnarray}\label{Functg}
\mbox{argmin}_{x \in C} g_D(x) = \mbox{sol}(VI(F, C)),
\end{eqnarray}
\noindent where sol($VI(F, C)$) is the solution set of the
variational inequality. Thus using the properties of the gap
function we can also write
\begin{eqnarray*}
\mbox{sol}(VI(F, C))= \{ x \in \mathbb{R}^n : g_D (x) = 0 , \  x \in C
\} .
\end{eqnarray*}
Hence, when the objective function $ f $ is convex,  we can write the
problem (SMPEC) equivalently as the following non-smooth simple bilevel programming problem (SBP1),
\begin{eqnarray*}
\min f(x) \quad\mbox{subject to}\quad \mbox{argmin}_{x \in C} g_D(x).
\end{eqnarray*}
Since $ \inf\limits_C g_D =0 $, the problem (SBP1) can be re-written as the following single level convex optimization problem (r-SMPEC)
\begin{eqnarray} \label{rSMPEC}
\min f(x) \quad\mbox{subject to}\quad g_D (x) \le 0, \quad x \in C.
\end{eqnarray}
\noindent However,  it is important to note that Slater constraint
qualification is never satisfied for the problem (r-SMPEC). Using
the reformulation (r-SMPEC) we have the following optimality
conditions. Before we write down the optimality conditions we would
like to note that $ g_D $ is in general a non-differentiable convex
function and its subdifferential can be easily computed if $ C $ is
a compact convex set. This is achieved by applying Danskin's
formula \cite{Dan} which shows that
\begin{eqnarray*}
\partial g_D(x) = \mbox{conv} \{ F(y) : y \in Y(x)\},
\end{eqnarray*}
\noindent where $ Y(x) = \{ y \in C : g_D(x) = \langle F(y), x - y
\rangle \}= \mbox{argmax}_{y\in C} \langle F(y), x - y \rangle $.

\begin{Remark}
It is worth observing also that the function $g_D$ in (\ref{Functg})
is nothing else than the supremum of the family $\psi_y$ over $y \in
C$, i.e., $g_D(x) = \sup_{y \in C} \psi_y (x)$, and the equivalence
of the problem (\ref{32}) and (r-SMPEC) is well-known \cite{GL98}.
The most interesting feature in the previous reconstruction lies on
the fact that $g_D$ is exactly the dual gap function of $VI(F,C)$,
which ensures that $g_D(x) \geq 0$ for all $x \in C$ and $S= \{ x \in C\; :\; g_D(x) = 0
\}$.
\end{Remark}

\noindent {\it Optimality condition under the weak basic constraint
qualification.} Let us now consider the following approach to the
problem (r-SMPEC). We shall first write down the standard optimality
condition in terms of the normal cone associated with the problem
(r-SMPEC)  and use the weak basic constraint qualification due to
Dempe and Zemkoho \cite{DeZe11} to write down an estimation for the
normal cone. In fact it has been shown in \cite{DeZe11} through
examples that for the reformulated version of a simple bilevel
problem even though the Slater condition fails the weak basic
constraint qualification can still hold. In fact if the Slater
condition hold true in the case of (r-SMPEC) then it could be
equivalently written as
\begin{eqnarray*}
\partial g_D (\bar{x}) \cap (-N_C(\bar{x})) = \emptyset.
\end{eqnarray*}
\noindent This formulation is what is known as basic constraint
qualification in the literature when we have just one functional
constraint. See Rockafellar and Wets \cite{rocw} for more details on
this. Since the Slater condition does not hold of course the basic
constraint qualification does not hold true for (r-SMPEC). However a
slightly weaker version of the basic constraint qualification as
given by \cite{DeZe11} is as follows. We will say that the problem
(r-SMPEC) satisfies the {\it weak basic constraint qualification} (in brief,  {\it weak BCQ})  at
$\bar{x}$ if
\begin{eqnarray*}
\partial g_D(\bar{x}) \cap (-\bd N_C(\bar{x})) = \emptyset,
\end{eqnarray*}
where $\bd A$ denotes the (relative) boundary of the set $A$.
\noindent Note that if the basic constraint qualification holds then
the above relation also holds and this justifies the name weak basic
constraint qualification. Note that what it implies is that if the
basic constraint qualification does not hold but the weak basic
constraint qualification holds then the subdifferential $ \partial
g_D(\bar{x}) $ must belong to the relative interior or without loss of
generality the interior of the negative of the normal cone to $ C $
at $\bar{x} $. We now present the following result.
\begin{Theorem}
Consider the problem (r-SMPEC) where $C$ is a compact, convex subset
of $\R^n$.  Let  $ \bar{x} $ be an optimal solution. Assume that the
weak  BCQ holds at $\bar{x} $. Then there exist scalars $
\beta_i \ge 0 $, $ i = 1, \ldots, (n +1) $ and vectors $ y_i $, with
$ y_i \in Y(\bar{x}) $, $ i = 1, \ldots ,(n+1)$ such that
\begin{eqnarray*}
- \sum_{i = 1}^{n +1} \beta_i F( y_i) \in \partial f(\bar{x})+
N_C(\bar{x}).
\end{eqnarray*}
\noindent Further,  for any $ \bar{x} \in C $ if there exist scalars
$ \beta_i \ge 0 $,   and vectors $ y_i $, with $ y_i \in Y(\bar{x})
$, $ i = 1, \ldots ,(n+1)$ such that the above inclusion holds, then
$\bar{x} $ is an optimal solution of (r-SMPEC).
\end{Theorem}

\noindent{\bf Proof.} If $ \bar{x} $ is a solution of the problem
(r-SMPEC) then from the standard optimality conditions for a convex
optimization problem we have that
\begin{eqnarray*}
0 \in \partial f(\bar{x}) + N_S(\bar{x}).
\end{eqnarray*}
\noindent Since the weak BCQ holds at $\bar{x}$, using Lemma 3.3
in Dempe and Zemkoho \cite{DeZe11}, we have that
\begin{eqnarray*}
N_S(\bar{x}) \subset \bigcup_{\alpha \ge 0} \alpha \partial
g_D(\bar{x}) + N_C(\bar{x}).
\end{eqnarray*}
\noindent Thus there exists $\lambda \ge 0 $ such that
\begin{eqnarray*}
0 \in \partial f(\bar{x}) + \lambda \partial g_D(\bar{x}) +
N_C(\bar{x}).
\end{eqnarray*}
\noindent Since $ C $ is compact we know that \cite{Dan}
\begin{eqnarray*}
\partial g_D(\bar{x}) = \mbox{conv} \{
F(y) : y \in Y(\bar{x}) \}.
\end{eqnarray*}
\noindent By Carath\'{e}odory's theorem,  there exist
scalars $ \mu_i \ge 0 $, $i = 1, \ldots , (n +1) $ with $
\sum_{i=1}^{n+1} \mu_i = 1 $ and vectors  $ y_i $, $ i= 1, \ldots,
(n + 1) $ with each $ y_i \in Y(\bar{x}) $   such that
\begin{eqnarray*}
- \lambda \sum_{i = 1}^{n +1} \mu_i F( y_i) \in \partial f(\bar{x})+
N_C(\bar{x}).
\end{eqnarray*}
\noindent The result is established by setting $ \beta_i = \lambda
\mu_i $ for each $ i = 1, \ldots, (n + 1) $.\\
\noindent The proof of the converse is omitted since it is straightforward. \hfill $\Box$

\begin{Example} \rm
This is an example where weak BCQ does not hold.
Let $F(x)=x$ and $C=[-1,1]$. Then $0$ is a solution of the $VI(F,C)$.
We consider the dual gap function $g_D$,
\begin{eqnarray*}
g_D(x)= \sup\limits_{y \in C} \langle y, x-y \rangle .
\end{eqnarray*}
Then
\begin{eqnarray*}
\partial g_D(0)=\{0\}.
\end{eqnarray*}
Also,  we have that $0 \in \bd N_C(0)$.
Hence,  $\partial g_D (0) \bigcap (-\bd N_C(0))=\{0\} \neq \emptyset$. \hfill $\Box$
\end{Example}

\begin{Example} \rm
This is an example where weak BCQ holds, modified from the example given by Dempe and Zemkoho \cite{DeZe11}.

We consider the simple convex bilevel programming problem.
\begin{eqnarray*}
&&\min x^2+y^2\\
\mbox{such that} && (x,y) \in S:= \mbox{argmin} \{x+y \, :\,  (x,y)\in \Omega\}.
\end{eqnarray*}
We have $f(x,y):= x+y$ and $\Omega :=\{ (x,y)\in \mathbb{R}^2 \, :\,  0 \leq x \leq 1 , 0 \leq y \leq 1 \}$. The point $(x_1,y_1)=(0,0)$ is the unique optimal solution of the problem and $N_{\Omega}(0,0)= \mathbb{R}_-^2$.
As solving the lower level optimization problem is the same as solving  $VI(F, \Omega)$, where $F(x,y)=\nabla f (x,y)  $  and  $\nabla f (x,y) = (1,1)$; the corresponding dual gap function is
\begin{eqnarray*}
g_D(x_1,y_1) = \sup \limits _{(x,y) \in \Omega} \langle (1,1), (x_1,y_1)- (x,y) \rangle .
\end{eqnarray*}
Then
\begin{eqnarray*}
g_D(0,0) = \sup \limits _{(x,y) \in \Omega} [-x-y] =0 .
\end{eqnarray*}
Now by Danskin's theorem \cite{Dan} we get that $\partial g_D(0,0)=\{(1,1)\}$.\\
Further,
\begin{eqnarray*}
(1,1) \notin \{(x,0) \, :\,  x \geq0\} \cup \{(0,y) \, :\,  y \geq 0\} = -\bd N_{\Omega}(0,0).
\end{eqnarray*}
Clearly, $\partial g_D(0,0) \cap (-\bd N_{\Omega}(0,0))= \emptyset $. \hfill $\Box$
\end{Example}

\subsection{Calmness and Optimality}
The notion of calmness is a very fundamental notion used for
deriving Lagrange multiplier rules for various classes of
optimization problems, for example see Clarke \cite{cla83}. In this
section we shall show that the notion of calmness used by Henrion,
Jourani and Outrata \cite{rene2002} can be applied to derive a
necessary and sufficient condition for the problem (SMPEC). Let us
first begin by defining the notion of calmness as given in  \cite{rene2002}. Let $ M : \mathbb{R}
\rightrightarrows \mathbb{R}^n $ be a set-valued map and let $(
\bar{y}, \bar{x}) \in \gph M $, where $ \gph M $ denotes the graph of
the set-valued mapping $ M $. Then $ M $ is said to be calm at $ (
\bar{y}, \bar{x}) $ if there exist neighborhoods $ V $ of $ \bar{y}$
and $ U $ of $ \bar{x} $ and a real number $ L > 0 $ such that
\begin{eqnarray*}
d( x, M(\bar{y})) \le L d(y, \bar{y}) \quad \forall y \in V, \quad
\forall x \in M(y) \cap U.
\end{eqnarray*}
Now consider the problem (r-SMPEC), where $ f $ is  convex and $ F
$ is a monotone map. It is not very difficult to observe that
(r-SMPEC) can be embedded in a larger family of problems $P(y)$, where
$ y \in \mathbb{R} $,  given by
\begin{eqnarray*}
\min f(x) \quad\mbox{subject to}\quad x \in M(y),
\end{eqnarray*}
\noindent where
\begin{equation} \label{peq1}
 M(y) = \{ x \in C : g_D(x) + y \le 0 \} .
\end{equation}
\noindent Note that $ P(0) $ is nothing but (r-SMPEC). If $ \bar{x}
$ is a solution of (r-SMPEC) then from the standard optimality
conditions  in convex optimization we have
\begin{eqnarray*}
0 \in \partial f(\bar{x}) + N_{M(0)}( \bar{x}).
\end{eqnarray*}
\noindent Let us now assume that the set valued map $ M : \mathbb{R}
\rightrightarrows \mathbb{R}^n $ given by (\ref{peq1}) is calm at $(0,
\bar{x}) \in \gph M $. Then by using Theorem 4.1 in \cite{rene2002}
we conclude that
\begin{eqnarray*}
N_{M(0)} (\bar{x}) \subseteq \bigcup_{y^* \geq 0} D^*
g_D(\bar{x})(y^*) + N_C(\bar{x}),
\end{eqnarray*}
where $D^*$ refers to Mordukhovich’s coderivative, see \cite[Definition 1.32]{boris1}.
\noindent It now follows from \cite{boris1}  that
\begin{eqnarray*}
 D^*g_D(\bar{x})(y^*) = y^* \partial g_D(\bar{x}),
 \end{eqnarray*}
and thus,
\begin{eqnarray*}
0 \in \partial f(\bar{x}) + y^* \partial g_D(\bar{x}) +
N_C(\bar{x}).
\end{eqnarray*}
It is now simple to see that if there is an $ \bar{x} \in S $ and $
y^* \ge 0 $ such that the above inclusion holds then $ \bar{x} $
solves (SMPEC). We sum up the discussion above in the following
theorem.
\begin{Theorem}
Consider the problem (r-SMPEC) where $ f $ is convex and $ F $ is
continuous and monotone. Let $ \bar{x} $ be feasible to (r-SMPEC).
Let us consider the set valued map $ M : \mathbb{R} \rightrightarrows
\mathbb{R}^n $ given by (\ref{peq1}). Assume that $ M $ is calm at $
(0, \bar{x}) \in gph M $. Then $ \bar{x} $ is a minimizer of
(r-SMPEC) if and only there exists $ y^* \ge 0 $ such that
\begin{eqnarray*}
0 \in \partial f(\bar{x}) + y^* \partial g_D(\bar{x}) +
N_C(\bar{x}).
\end{eqnarray*}
\end{Theorem}
\noindent In fact if $ F $ is continuous 
 and C is compact then we have \[ \partial g_D(\bar{x}) = \mbox{conv}_{y \in Y(\bar{x})}
\{F(y)\},\textrm{ where } Y(\bar{x}) = \mbox{argmax}_C \langle F(y),
\bar{x} - y \rangle .\]
\section{Optimality conditions and schematic algorithm for \\(SMPEC)}
We divide this section into two subsections. In the first subsection we study sequential optimality conditions, for which no constraint qualifications are needed. In the second subsection we present a schematic algorithm using the dual gap function and show how the sequential optimality conditions play a role in the convergence analysis.
\subsection{Sequential optimality conditions}

We shall now present sequential optimality conditions that hold
without any constraint qualification nor qualification conditions
for (SMPEC). We shall present two approaches. The first one concerns
conjugate theory of convex analysis and approximate subdifferentials
while for the second one,  a sequential chain rule in \cite{tl} is
used to derive the sequential optimality condition for the case
where the convex constrained set $C$ is explicitly defined by a
 convex inequality. We start with the first approach.

\medskip
{\it The case where the constraint set C is not explicitly given.}  Let $A$ be the
feasible set of the problem (r-SMPEC). It is worth noting that by
the formulation of the (r-SMPEC), if $x \in A$ then $g_D(x) = 0$. We
first need the following lemma.

\begin{Lemma} \label{lem51} Assume that $C$ is a nonvoid compact, convex subset, $A \not=\emptyset$ and $\bar x \in A$. If $v \in \partial \delta_A (\bar x)$
then for each $k \in \N$, there exist $\lambda_k \in \R_+$, $x_k \in
\dom  g_D$, $y_k \in C$, $v_k \in
\partial (\lambda_k g_D)(x_k)$, $w_k \in N_C(y_k)$ such that as $k
\rightarrow +\infty$, it holds
\begin{eqnarray*}
&&v_k + w_k \rightarrow v, \\
&& x_k \rar \bar x , \ y_k \rar \bar x,  \\
&& \lambda_k g_D (x_k) - \la v_k, x_k - \bar x \ra \rar 0,  \ {\rm
and} \ \la  w_k, y_k - \bar x \ra \rar 0.
\end{eqnarray*}
\end{Lemma}

\noindent{ \bf Proof.}  We first note that if $v \in \partial
\delta_A (\bar x)$ then by (\ref{epifstar}), $(v, \la v,\bar x\ra ) \in
\epi \delta_A^\ast$. From \cite[Lemma 2.1]{jld}, we have
\[
\epi \delta_A^\ast = \mathrm {cl}\big( \bigcup\limits _{\lambda \in
\R_{+}}\mbox{\rm epi } (\lambda g_D)^{\ast }+\mbox{\rm epi }\delta
_{C}^{\ast }\big),
\]
and thus, for each $k\in \N$, there exist $\lambda_k \in \R_+$,
$(\bar v_k, \alpha_k) \in \epi (\lambda_k g_D)^{\ast }$, $(\bar w_k,
\beta_k) \in \epi \delta_C^\ast$ such that
$
\bar v_k + \bar w_k \rar v, \ \alpha_k+ \beta_k \rar \la v,\bar x\ra .
$
Using (\ref{epifstar}) one more time, applied to the functions
$(\lambda_k g_D)^\ast$ and $\delta_C^\ast$, for each $k$, there
exist $\eta_k, \xi_k  \in \R_+$ such that
\begin{eqnarray*}
&&\bar v_k \in \partial_{\eta_k} (\lambda_k g_D)(\bar x),  \  \ \alpha_k = \la \bar v_k,  \bar x \ra + \eta_k - \lambda_k
g_D(\bar x) = \la \bar v_k,  \bar x \ra + \eta_k\\
&& \bar w_k \in \partial_{\xi_k} \delta_C(\bar x),\  \mathrm{ and } \
\beta_k = \la \bar w_k, \bar x \ra + \xi_k.
\end{eqnarray*}
Since $\bar v_k + \bar w_k \rar v, \ \alpha_k + \beta_k \rar \langle v,\bar
x\rangle$, the previous inequalities imply that $\eta_k \downarrow 0$,
$\xi_k \downarrow 0$, and $\lambda_k g_D(\bar x) \rar 0$ as $k \rar + \infty$.

We are now ready to apply Brondsted-Rockafellar theorem (Proposition
\ref{proposition2.1}) to get the desired result. Indeed, it follows
from Brondsted-Rockafellar theorem, applied to $\partial_{\eta_k}
(\lambda_k g_D)(\bar x)$ and $\partial_{\xi_k} \delta_C(\bar x)$, $k
\in \N$, that there exist $x_k \in \dom(\lambda_k g_D)$, $y_k \in
C$, $v_k \in \partial (\lambda_k g_D)(x_k)$, $w_k \in N_C(y_k)$
satisfying
\begin{eqnarray*}
&&\Vert x_k - \bar x\Vert \leq \sqrt{\eta_k} ,\ \  \Vert v_k - \bar v_k\Vert \leq \sqrt{\eta_k},  \ \
\vert \lambda_k g_D(x_k) - \la v_k, x_k - \bar x\ra -\lambda_k g_D(\bar x) \vert \leq 2 \eta_k , \\
&&\Vert y_k - \bar x \Vert \leq \sqrt{\xi_k}, \ \
 \Vert w_k - \bar w_k\Vert \leq \sqrt{\xi_k} , \ \ \vert \la w_k, y_k - \bar x \ra \leq 2 \xi_k.
\end{eqnarray*}
Taking into account the fact that $\bar v_k + \bar w_k \rar v$, $\eta_k \rar 0$, $\xi_k \rar 0$, $\lambda_k g_D(\bar x) \rar 0$ as $k \rar + \infty$, we get
(as $k \rar + \infty$)
\begin{eqnarray*}
&& x_k \rar \bar x, \ y_k \rar \bar x ,  \ \
 v_k + w_k \rar v, \\
&& (\lambda_k g_D)(x_k) - \la v_k, x_k - \bar x\ra \rar 0,\\
&&\la w_k, y_k - \bar x\ra \rar 0,
\end{eqnarray*}
which is desired.  \hfill $\Box$

\begin{Theorem} \label{seqopt}
Let us consider the problem (SMPEC) where $ C $ is a  compact, convex subset, $ F
$ is monotone and continuous, and $ f $ be a finite-valued convex
function. Then $ \bar{x} $ is a solution of (SMPEC) if and only if
there exists $u \in \partial f(\bar x)$ and, for each $k \in \N$,
there exist
 $\lambda_k \in \R_+$, $x_k \in \dom (\lambda_k g_D)$, $y_k \in C$,
 $w_k \in N_C(y_k)$, $ \mu^k_i \ge 0 $ and vectors $ y^k_i \in Y(x_k) $  $ i = 1,
2, \ldots n +1 $ with $ \sum_{i  = 1}^{n + 1} \mu^k_i = 1 $ such
that, as $k \rar +\infty$,
\begin{eqnarray*}
&& u + \lambda_k \sum_{i = 1}^{n + 1} \mu^k_i F(y^k_i) + w_k \rar 0, \\
&& x_k \rar \bar x, \ y_k \rar \bar x ,  \\
&& (\lambda_k g_D)(x_k) - \la \lambda_k \sum_{i = 1}^{n + 1} \mu^k_i F(y^k_i), x_k - \bar x\ra \rar 0,\\
&&\la w_k, y_k - \bar x\ra \rar 0 .
\end{eqnarray*}

\end{Theorem}

\noindent{ \bf Proof.} The problem (r-SMPEC) can be reformulated as
an unconstrained convex problem (P):
\[
\inf\limits_{x \in \R^n} (f(x) + \delta_A(x)),
\]
where $A$ is the feasible set of (r-SMPEC), i.e., $A = C \cap
g_D^{-1} (0)$.

$\bullet$ {\it Necessary condition}. Assume that  $\bar x$ is a
solution of (r-SMPEC). Then it is a solution of (P).  Since $f$ is
continuous, it is continuous at $\bar x \in A$, and hence, by a
standard result in convex analysis, one has,
\[
0 \in \partial f(\bar x) + \partial \delta_A(\bar x),
\]
which ensure that there exist $u \in \partial f(\bar x)$, $v \in
\partial \delta_A(\bar x)$ such that $ u + v = 0$.  Now, since $v \in
\partial \delta_A(\bar x)$, by  Lemma \ref{lem51},
for each $k \in \N$, there exist $\lambda_k \in \R_+$, $x_k \in \dom
g_D$, $y_k \in C$, $v_k \in
\partial (\lambda_k g_D)(x_k)$, $w_k \in N_C(y_k)$ such that as $k
\rightarrow +\infty$, it holds,
\begin{eqnarray*}
&&u+ v_k + w_k \rightarrow 0, \\
&& x_k \rar \bar x , \ y_k \rar \bar x \\
&& \lambda_k g_D (x_k) - \la v_k, x_k - \bar x \ra \rar 0 \ {\rm
and} \ \la  w_k, y_k - \bar x \ra \rar 0.
\end{eqnarray*}
It is worth noting that
\begin{eqnarray*}
v_k \in \lambda_k \partial g_D(x_k) = \lambda_k \mbox{conv} \{ F(y)
: y \in Y(x_k) \}.
\end{eqnarray*}
 Using Carath\'{e}odory's theorem we can conclude that
there exists scalars $ \mu^k_i \ge 0 $ and vectors $ y^k_i \in
Y(x_k) $  $ i = 1, 2, \ldots n +1 $ with $ \sum_{i  = 1}^{n + 1}
\mu^k_i = 1 $ such that
\begin{eqnarray*}
v_k = \lambda_k \sum_{i = 1}^{n + 1} \mu^k_i F(y^k_i).
\end{eqnarray*}
This completes the necessary part.

$\bullet$ {\it Sufficient condition}. Assume the existence of $u \in
\partial f(\bar x)$ and the sequences satisfying the conditions of
the theorem. Let $x\in A$ be an arbitrary feasible point of
(r-SMPEC) (and hence,  $g_D(x) = 0$). Put further $v_k:= \lambda_k
\sum_{i = 1}^{n + 1} \mu^k_i F(y^k_i)$. Then $v_k \in \lambda_k
\partial g_D(x_k)$. From the fact that $u \in
\partial f(\bar x)$, $v_k \in
\lambda_k \partial g_D(x_k)$, $y_k \in C$, and
 $w_k \in N_C(y_k)$, we get
 \begin{eqnarray*}
f(x) - f(\bar x) &\geq& \la u, x - \bar x\ra \\
\lambda_k g_D(x) - \lambda_k g_D(x_k) &\geq& \la v_k, x - x_k\ra
\\
0 &\geq& \la w_k, x - y_k\ra ,
 \end{eqnarray*}
which entails
\begin{eqnarray*}\label{proof1}
f(x) - f(\bar x)  &\geq& \la u, x - \bar x\ra + \la v_k, x - x_k\ra
+ \la w_k, x - y_k\ra + \lambda_k g_D(x_k)\\
&=& \la u, x - \bar x\ra + \la v_k, x - \bar x\ra + \la w_k, x -
\bar x \ra  \\
&& +   \big( \lambda_k g_D(x_k) - \la v_k, x_k - \bar x\ra \big) -
\la w_k, y_k - \bar x\ra.
\end{eqnarray*}
Taking the limit when $k \rar + \infty$ the right hand side of the
previous inequality tends to zero by assumption. Consequently, $f(x)
- f(\bar x) \geq 0$ holds with arbitrary feasible point $x \in A$
which shows that $\bar x$ is a solution of (r-SMPEC), and hence, a
solution of (SMPEC). The proof is complete. \hfill $\Box$

\medskip

{\it The case with an explicit presentation of the  constrained set
$C$.}  In the optimality conditions that we present below we
consider  that set $ C $ is explicitly defined by a convex
inequality constraint, i.e.,
\begin{eqnarray*}
C = \{ x \in \mathbb{R}^n : h(x) \le 0 \},
\end{eqnarray*}
\noindent where $ h : \mathbb{R}^n \rightarrow \mathbb{R} $ is a
convex function. Let us note that this representation is quite
general since any convex set that is represented in terms of finite
number of convex inequality constraints can be represented in the
above form with the function $ h $ corresponding to the max function
of the associated convex functions.\\
\noindent Now we will present the sequential optimality conditions
using a sequential chain rule due to Thibault \cite{tl}.
\begin{Theorem}\cite{tl} \label{c-rule}
Consider the vector-valued function $ H : \mathbb{R}^n \rightarrow
\mathbb{R}^m $ whose each component is convex and let $ \phi :
\mathbb{R}^m \rightarrow \mathbb{R} \bigcup \{+ \infty\} $ be a proper,
lower-semicontinuous convex function which is increasing on $
H(\mathbb{R}^n) + \mathbb{R}^m_+ $. Then for $ H(\bar{x}) \in dom\,
\phi $, the vector $ \xi \in \partial( \phi \circ H)(\bar{x}) $ if
and only if there exist sequences $ \eta_k, y_k, \xi_k$ with $
\eta_k \in \partial \phi(y_k) $ and $ \xi_k \in \partial (\eta_k
H)(x_k) $ satisfying
\begin{enumerate}
\item[(i)]\ $ x_k \rightarrow \bar{x} $,
\item[(ii)] \ $ y_k \rightarrow H( \bar{x}) $,
\item[(iii)] \ $ \xi_k \rightarrow \xi $,
\item[(iv)]\ $ \phi(y_k) - \langle \eta_k, y_k -H(\bar{x})\rangle
\rightarrow \phi(H(\bar{x}))$,
\item[(v)] \ $ \langle \eta_k, H(x_k) - H(\bar{x})\rangle \rightarrow 0 $.
\end{enumerate}
\end{Theorem}
\noindent We would also like to note that for any $ \eta \in
\mathbb{R}^m $ the symbol $(\eta H)$ means the function which is
evaluated at each $ x $ as
$
(\eta H) (x) = \langle \eta, H(x) \rangle.
$

\begin{Theorem}
Let us consider the problem (SMPEC) where $ C= \{ x \in \mathbb{R}^n : h(x) \le 0 \} $ is compact, $ h : \mathbb{R}^n \rightarrow \mathbb{R} $ is a
convex function, $ F $
is monotone and continuous and $ f $ be a finite-valued convex
function. Then $ \bar{x} $ is a solution of (SMPEC) if and only if
there exist sequences $ x_k \rightarrow \bar{x} $, $ y_k = ( y^1_k,
y^2_k)\in - \mathbb{R}^2_+ $, $y^1_k  \rightarrow 0 $, $ y^2_k
\rightarrow h(\bar{x})$, $\lambda_k = (\lambda^1_k, \lambda^2_k) \in
\mathbb{R}^2_+ $, and $\xi \in\partial f(\bar{x}) $, and a sequence
$ v^k \in \partial h(x_k) $ and for every $ k \in \mathbb{N} $ there
exists scalars $ \mu^k_i \ge 0 $ and vectors $ y^k_i \in Y(x_k) $  $
i = 1, 2, \ldots n +1 $ with $ \sum_{i  = 1}^{n + 1} \mu^k_i = 1 $
such that,  as $ k \rightarrow \infty$,
\begin{enumerate}
\item[(i)]\  $ \xi + \lambda_k \sum_{i = 1}^{n + 1} \mu^k_i F(y^k_i) +
\lambda^2_k v^k \rightarrow 0$,
\item[(ii)]\ $\lambda^1_k y^1_k + \lambda^2_k y^2_k - \lambda^2_k h(\bar{x})
\rightarrow 0 $,
\item[(iii)]\ $ \lambda^1_k g_D(x^k) + \lambda^2_k (h(x_k) - h(\bar{x}))
\rightarrow 0 $.
\end{enumerate}
\end{Theorem}

\noindent{ \bf Proof. } Let us begin by defining the vector function
$ G : \mathbb{R}^n \rightarrow \mathbb{R}^2 $ given as
\begin{eqnarray*}
G(x) = (g_D(x), h(x)).
\end{eqnarray*}
\noindent Now if $\bar{x} $ solves (SMPEC) then $\bar{x} $ also
solves (r-SMPEC) which can now be written as
\begin{eqnarray*}
\min f(x), \quad\quad  G(x) \in -\mathbb{R}^2_+.
\end{eqnarray*}
\noindent Thus $\bar{x}$ is also a solution of the problem
\begin{eqnarray*}
\min_{x \in \mathbb{R}^n} (f + \delta_{-\mathbb{R}^2_+} \circ G)(x).
\end{eqnarray*}
\noindent Thus from standard optimality conditions in convex
optimization we have
\begin{eqnarray*}
0 \in \partial ( f + \delta_{-\mathbb{R}^2_+} \circ G)(\bar{x}).
\end{eqnarray*}
\noindent It is easy to see that $ \delta_{-\mathbb{R}^2_+} $ is
non-decreasing over the set $ G(\mathbb{R}^n) + \mathbb{R}^2_+ $.
Thus our first job is to show that $ \delta_{-\mathbb{R}^2_+} \circ
G $ is a convex proper function. Note that since both $ g_D $ and $
h $ are  convex functions we have for any $ x, y \in
\mathbb{R}^n $ and $ \lambda \in [0, 1] $
\begin{eqnarray*}
\lambda G(x) + ( 1 - \lambda) G(y) - G( \lambda x + (1 - \lambda) y)
\in \mathbb{R}^2_+.
\end{eqnarray*}
\noindent This clearly shows that
\begin{eqnarray*}
\lambda G(x) + (1 - \lambda) G(y) \in G(\mathbb{R}^n) +
\mathbb{R}^2_+.
\end{eqnarray*}
\noindent It is obvious  that $ G(\lambda x + ( 1 - \lambda)
y) \in G(\mathbb{R}^n) + \mathbb{R}^2_+ $. Thus, using that fact that
$ \delta_{-\mathbb{R}^2_+} $ is non-decreasing on $ G(\mathbb{R}^n)
+ \mathbb{R}^2_+ $,  we have
\begin{eqnarray*}
\delta_{-\mathbb{R}^2_+} \circ G (\lambda x + ( 1 -\lambda) y) \le
\lambda \delta_{-\mathbb{R}^2_+} \circ G (x) + ( 1 - \lambda)
\delta_{-\mathbb{R}^2_+} \circ G (y).
\end{eqnarray*}
\noindent Since the feasible set $ S $ which is the solution of $
VI(F, C) $ is assumed to be nonempty,  we see that
$\delta_{-\mathbb{R}^2_+} \circ G $ is a proper lower-semicontinuous
convex function. On the other hand, as  $ f $ is finite-valued,  it is continuous on
$\mathbb{R}^n $ and thus,  invoking the sum rule for subdifferentials
of a convex function,  we have
\begin{eqnarray*}
0 \in \partial f(\bar{x}) + \partial (\delta_{-\mathbb{R}^2_+} \circ
G)(\bar{x}),
\end{eqnarray*}
\noindent and so,  there exists $ \xi \in \partial f(\bar{x}) $ and $
\hat{\xi} \in \partial (\delta_{-\mathbb{R}^2_+} \circ G)(\bar{x}) $
such that
\begin{equation}\label{neweq10}
\xi + \hat{\xi} = 0.
\end{equation}
\noindent Now as $ G(\bar{x}) \in \dom \delta_{-\mathbb{R}^2_+} $,
using Theorem \ref{c-rule} we conclude that there exists sequences $
\lambda_k = (\lambda^1_k, \lambda^2_k) \in \mathbb{R}^2_+ $, $ y_k =
( y^1_k, y^2_k ) \in -\mathbb{R}^2_+ $, $ \xi_k $ and $ x_k $ with
$\lambda_k \in \partial \delta_{-\mathbb{R}^2_+} (y_k) =
N_{-\mathbb{R}^2_+} (y_k)$ and $ \xi_k \in \partial (\lambda_k
G)(x_k) $ such that

\begin{enumerate}
\item [(i)]\  $ x_k \rightarrow \bar{x} $

\item [(ii)]\  $ y^1_k \rightarrow 0$ and $ y^2_k \rightarrow h(\bar{x})$

\item [(iii)]\  $\xi_k \rightarrow \hat{\xi}$

\item [(iv)]\  $ \delta_{-\mathbb{R}^2_+} (y_k) - \langle \lambda_k, y_k
- G(\bar{x})\rangle \rightarrow
\delta_{-\mathbb{R}^2_+}(G(\bar{x}))$

\item[(v)]\  $\langle \lambda_k, G(x_k) - G(\bar{x})\rangle \rightarrow
0$.
\end{enumerate}

 \noindent It is clear that $ y_k \in
-\mathbb{R}^2_+ $ and further (ii), (iv) and (v) above reduce
to
\begin{eqnarray*}
&&\lambda^1_k g_D(x_k) + \lambda^2_k( h(x_k) - h(\bar{x})) \rightarrow
0    \ \ \text{and} \\
&&\lambda^1_k y^1_k + \lambda^2_k y^2_k - \lambda^2_k
h(\bar{x})\rightarrow 0.
\end{eqnarray*}
\noindent Now noting that $ \lambda_k \in N_{-\mathbb{R}^2_+} (y_k)$
we have
\begin{eqnarray*}
\langle \lambda_k, y - y_k \rangle \le 0 \quad \forall y \in
-\mathbb{R}^2_+.
\end{eqnarray*}
Putting $ y = 0 $ and $ y = 2 y_k $ in the above expression for
the normal cone we see that $ \langle \lambda_k, y_k \rangle = 0 $.
This immediately shows that $ \langle \lambda_k , y \rangle \le 0 $
for all $ y \in -\mathbb{R}^2_+ $. Hence,  $
\lambda_k \in \mathbb{R}^2_+ $. Now again using the sum rule for the
subdifferential of convex functions we see that
\begin{eqnarray*}
\xi_k \in \lambda^1_k \partial g_D (x_k) + \lambda^2_k \partial
h(x_k) .
\end{eqnarray*}
From (\ref{neweq10}) we conclude that there exists $ v^k \in
\partial h(x_k) $ and for every $ k \in \mathbb{N} $ there
exists scalars $ \mu^k_i \ge 0 $ and vectors $ y^k_i \in Y(x_k) $,   $
i = 1, 2, \ldots n +1 $ with $ \sum_{i  = 1}^{n + 1} \mu^k_i = 1 $
such that (as $ k \rightarrow \infty$)
\begin{eqnarray*}
 \xi + \lambda^1_k \sum_{i = 1}^{n + 1} \mu^k_i F(y^k_i) + \lambda^2_k v^k
\rightarrow 0.
\end{eqnarray*}
\noindent This completes the necessary part. For the converse,  it is
simple to see that from the given conclusion we can immediately
derive that
\begin{eqnarray*}
-\xi \in \partial (\delta_{-\mathbb{R}^2_+} \circ G)(\bar{x}).
\end{eqnarray*}
\noindent This of course shows that
\begin{eqnarray*}
0 \in \partial f(\bar{x}) + \partial (\delta_{-\mathbb{R}^2_+} \circ
G)(\bar{x}),
\end{eqnarray*}
\noindent showing  that $\bar{x}$ solves (r-SMPEC) and hence, solves
(SMPEC). \hfill $\Box$

\subsection{ A schematic algorithm for the SMPEC problem}
In this section we will try to connect the sequential optimality conditions  with the numerical scheme we are proposing for the (SMPEC) problems. Further it was emphasized that we seek some connections between the sequential optimality conditions which are free of constraint qualifications and the numerical scheme we developed here for the (SMPEC) problem in \cite{part2}. The numerical scheme developed for the (SMPEC) problem in \cite{part2} uses quite different tools compared to ones that appear in the sequential optimality conditions. Thus apparently a connection is hard to find. We thus asked ourselves the following question. Is it possible at least from a theoretical perspective to develop a numerical scheme where sequential optimality conditions will play a central role. In fact we will discuss below a numerical scheme when the set $ C $ is convex and compact. We will see that Theorem \ref{seqopt} will play a pivotal role in it. Note that the (SMPEC) problem can be compactly written as
\begin{eqnarray*}
\min f(x), \quad \mbox{subject to} \quad x \in \mbox{sol}(VI(F,C))
\end{eqnarray*}
In fact if  sol$(VI(F,C)) \neq \emptyset $, we recall from section 3.2 that the (SMPEC) problem equivalently can be written as the problem (SBP1) given as,
\begin{eqnarray*}
\min f(x), \quad \mbox{ subject to } \quad x \in \mbox{argmin}_{x \in C} g_D(x)
\end{eqnarray*}
For the practical purposes it is imporatnt to have $ C $ to be compact set. Then it guarantees that $ g_D $ is finite.\\
\noindent  Here we present a schematic algorithm for the (SMPEC) problem where any accumulation point of the iterates is a solution of the (SMPEC) problem. This will be achieved by sequentially solving the problem (SBP1) by a kind of "{regularization}" approach. In the convergence analysis it is not easy to show directly that the iterates generated by the algorithm converges to a solution. In fact what we show here is that any limit point of the iterates satisfy the sequential optimality conditions. The algorithm is presented as follows:
\begin{tcolorbox}
\begin{description}
  \item[Step1] Let us choose $ x_0 \in C $, $ k = 0 $ and the optimality threshold $ \mu > 0$.
  \item[Step2] For any $ k \in \mathbb{N} $, $ x_k $ is a solution of the following optimization problem, let us call it as $ (P_k) $:
\begin{eqnarray*}
\min\limits_{x \in C } f(x) + \varepsilon_k g_D(x)
\end{eqnarray*}
where $ \varepsilon_k > 0 $ for all $ k \in \mathbb{N} $, $ \varepsilon_k \downarrow 0 $ and $ g_D(x) = \sup\limits_{y \in C} \langle F(y), x-y \rangle \geq 0 $.
  \item[Step3](Stopping Criteria) Choose $ x_k \in C $ as an approximate solution of (SMPEC) if $ g_D(x_k) < \mu $.
\end{description}
\end{tcolorbox}
\noindent  The approach that we seek here is not a penalization approach but is rather similar to a regularization scheme in optimization. The stopping criterion mentioned in Step 3 of the algorithm is motivated from the simple fact that if we get any iteration point $ x_k $ such that $ g_D(x_k) =0 $, then $ x_k $ will be a solution of the (SMPEC) problem. It is obtained as follows. Let us take any point $ \tilde{x} \in Sol(VI(F,C))$, then $ g_D(\tilde{x}) = 0 $ and also $ \tilde{x} \in C $. As $ x_k \in C $ is a solution of the problem $(P_k)$, we have
\begin{eqnarray*}
f(x_k) + \varepsilon_k g_D(x_k) \leq f(\tilde{x}) + \varepsilon_k g_D(\tilde{x})
\end{eqnarray*}
Using the fact that $ g_D(x_k) \geq 0 $ and $ g_D(\tilde{x}) = 0 $ we get that $ f(x_k) \leq f(\tilde{x})$. This is true for any $ \tilde{x} \in Sol(VI(F,C))$. Hence $ x_k $ is a solution of the (SMPEC) problem, if $ g_D(x_k) = 0 $. Thus from the practical point of view we can stop the algorithm if $ g_D (x_k) < \mu $, where $ \mu > 0 $ is a pre-set optimality threshold. 

\begin{Theorem}
Let us consider the (SMPEC) problem with $ C $ a compact, convex set. If $ \{x_k\} $ is the sequence generated by the algorithm mentioned above, then any accumulation point $ \bar{x} $ of $ \{x_k\} $ is a solution of the (SMPEC) problem.
\end{Theorem}
\noindent  {\bf Proof.}
Since $ x_k $ is a solution of the problem $(P_k)$ and $ C $ is compact, by using the famous Rockafellar-Pschenichnyi condition we obtain the following
\begin{eqnarray}
0 \in \partial f(x_k) + \varepsilon_k \partial g_D(x_k) + N_C(x_k).
\end{eqnarray}
Hence there exists $ u_k \in \partial f(x_k) $, $ v_k \in \partial g_D(x_k) $ and $ w_k \in N_C(x_k) $ such that
\begin{eqnarray}\label{op}
u_k + \varepsilon_k v_k + w_k =0.
\end{eqnarray}
Since $ v_k \in \partial g_D(x_k) $, using the subdifferential formula for $ \partial g_D(x_k) $ and the Caratheodory's theorem, there exists $ \mu_i^k \geq 0 $, $ \sum\limits_{i=1}^{n+1} \mu_i^k =1 $ and $ y_i^k \in Y(x_k) $ such that $ v_k = \sum\limits_{i=1}^{n+1} \mu_i^k y_i^k $. This implies that
\begin{eqnarray}\label{op1}
u_k + \varepsilon_k \sum_{i=1}^{n+1} \mu_i^k y_i^k + w_k =0.
\end{eqnarray}
As $ x_k \in C $ and compactness of $ C $ implies that $ \{ x_k \}  $ has a convergent subsequence. Without loss of generality, let us denote that subsequence by $ \{ x_k \} $ which converges to $ \bar{x} \in C $. Then using the locally boundedness property of the subdifferential we get that $ \{ v_k \} $ is a bounded sequence. As $ k \rightarrow \infty $, $ \varepsilon_k \rightarrow 0 $ and hence
\begin{eqnarray}\label{neg}
\lim\limits_{ k \rightarrow \infty }[ \varepsilon_k g_D(x_k) - \varepsilon_k \langle v_k, x_k - \bar{x} \rangle ] = 0.
\end{eqnarray}
Now, from \eqref{op} we get that $ \varepsilon_k v_k + w_k = -u_k $. Then
\begin{eqnarray*}
\varepsilon_k g_D(x_k) - \varepsilon_k \langle v_k, x_k - \bar{x} \rangle + \langle w_k , \bar{x} - x_k \rangle
&&= \varepsilon_k g_D(x_k) + \langle -u_k, \bar{x} - x_k \rangle\\
&& \geq \varepsilon_k g_D(x_k) + f(x_k) - f(\bar{x})
\end{eqnarray*}
as, convexity of $ f $ implies that $ f(\bar{x}) \geq f(x_k) + \langle u_k , \bar{x} - x_k \rangle $. Using the continuity of $ f $ and finite valuedness of $ g_D $ on $ C $, we get that
\begin{eqnarray*}
\lim\limits_{k \rightarrow \infty } [\varepsilon_k g_D(x_k) - \varepsilon_k \langle v_k, x_k - \bar{x} \rangle + \langle w_k , \bar{x} - x_k \rangle] \geq 0^.
\end{eqnarray*}
As $ w_k \in N_C(x_k)$, we have $ \langle w_k, \bar{x} - x_k \rangle \leq 0 $. This together with (\ref{neg}) implies that
\begin{eqnarray*}
\lim\limits_{ k \rightarrow \infty } [ \varepsilon_k g_D(x_k) - \varepsilon_k \langle v_k, x_k - \bar{x} \rangle + \langle w_k , \bar{x} - x_k \rangle ] \leq 0.
\end{eqnarray*}
Hence,
\begin{eqnarray}\label{eq}
\lim\limits_{k \rightarrow \infty } [\varepsilon_k g_D(x_k) - \varepsilon_k \langle v_k, x_k - \bar{x} \rangle + \langle w_k , \bar{x} - x_k \rangle] = 0.
\end{eqnarray}
Then from (\ref{neg}) and (\ref{eq}) we can conclude that as $ k \rightarrow \infty $
\begin{eqnarray}\label{op3}
 \langle w_k , \bar{x} - x_k \rangle \rightarrow 0.
\end{eqnarray}
Using the representation of $ v_k $, from (\ref{neg}) we get that
\begin{eqnarray}\label{op2}
\varepsilon_k g_D(x_k) - \langle \varepsilon_k \sum_{i=1}^{n+1} \mu_i^k y_i^k, x_k - \bar{x} \rangle \rightarrow 0.
\end{eqnarray}
Again using the locally boundedness property of subdifferential we know that $ \{ u_k \}$ is a bounded sequence. Let $ \{ u_k \}$ be a subsequence again denoted by  $ \{ u_k \}$ such that $ \lim\limits_{k \rightarrow \infty } u_k = u $, then $ u \in \partial f ( \bar{x}) $. Then from (\ref{op1}) we can conclude that as $ k \rightarrow \infty $
\begin{eqnarray}\label{op4}
u + \varepsilon_k \sum_{i=1}^{n+1} \mu_i^k y_i^k + w_k  \rightarrow 0.
\end{eqnarray}
Therefore $ x_k \rightarrow \bar{x} $ and (\ref{op4}), (\ref{op2}), (\ref{op3}) all together satisfy the optimality conditions mentioned in Theorem 4.1 (with $ y_k = x_k $ and $ \lambda_k = \varepsilon_k $). Hence $ \bar{x} $ is a solution of the SMPEC problem. This is true for any accumulation point of the sequence $ \{x_k\}$. Therefore we can conclude that all the accumulation points of the sequence generated by the algorithm are solutions of the (SMPEC) problem. \hfill $\Box$
\begin{Remark}
It is to be noted that the above scheme is a conceptual one. One of the major practical issues here is the computation of $ g_D $. In fact from a theoretical point of view we will just analyse how $ g_D $ can be replaced by a more computable gap function, in our future research on this class of problems. Another reason to call this algorithm schematic is we tacitly assume here that values of $ g_D $ and the elements of $ \partial g_D $ are computed through a black-box.  
\end{Remark}

\section{Conclusion} 
In this article we have tried to use several different approaches to develop necessary and sufficient optimality condition. A natural problem that we need to consider in our future research is the (MPGVI) problem, {\it i.e.} the (SMPEC) problem where the lower level has a GVI (Generalised VI), instead of VI defined by a maximal monotone map. Is it possible to develop necessary and sufficient optimality conditions for the (MPGVI)?  The very nature of GVI tells us that this may not achievable as a straight forward generalisation of the optimality conditions presented in this paper. We might also try to seek an answer as to why we need to devise necessary and sufficient optimality conditions through different approaches. One simple answer to this is that these conditions can be used as a benchmark to devise a stopping criteria for algorithms one might use to solve (SBP) and (SMPEC) problems. This issue is discussed in Part-II \cite{Part2}. The use of sequential optimality condition in the convergence analysis of a numerical scheme is demonstrated in section 4.2.\\
One of the key difficulties in analysing the (SMPEC) or (SBP) problem from the view of optimality conditions is the failure of the Slater condition even if we pose them equivalently as single level problems. In this article we had tried to overcome this by trying to use various weaker conditions. It might be interesting to see through examples if weaker constraints due to Abadie and Guignard hold for (SBP) problems. Thus even from the theoretical aspect of optimality conditions the (SBP) and (SMPEC) problems remain to be challenging problems to study.

\end{document}